\documentclass[12pt]{article}
\usepackage{graphicx} 
\usepackage{epstopdf}
\usepackage{subfigure}
\usepackage{color}
\usepackage[english]{babel}
\usepackage{amsmath,amsthm}
\usepackage{amsfonts}
\usepackage{booktabs}
\usepackage{amssymb}

\usepackage{graphicx}
\usepackage{subfigure}
\usepackage{tabularx}
\title{A note on the Nearly Dispersability of Odd Toroidal Grids}

\author { Xiaoxiang Yu, Zeling Shao, Zhiguo Li{$^*$}\\
{\small School of Science, Hebei University of Technology, Tianjin 300401, China}
\date{}
\footnote{Corresponding author. E-mail: zhiguolee@hebut.edu.cn}
}

\begin{document}
\baselineskip 0.65cm

\maketitle

\begin{abstract}
  The \emph{matching book thickness} $mbt(G)$ of $G$ is the minimum integer $m$ such that an $m$-page matching book embedding exists. A graph $G$ is called \emph{dispersable} if $mbt(G)=\Delta(G)$, \emph{nearly dispersable} if $mbt(G)=\Delta(G)+1$. Recently, the authors determined the nearly dispersability of odd toroidal grids $T_{s,t}$. In this note, we further present a brief proof for this result.

\bigskip
\noindent\textbf{Keywords:} Book embedding; Dispersability;  Toroidal grids; Layout

\noindent\textbf{2020 MSC.} 05C10
\end{abstract}

\section{Introduction}

The development of book embedding has been inspired and guided by computer science. At the same time, book embedding research also plays an increasingly important role in some directions of computer science, such as sorting with parallel stacks, fault-tolerant computing, and so on $[1]$. 
Let $\psi$ be a permutation of all vertices of a graph $G$, a \emph{layout} $\Psi=(G,\psi)$ for $G$ is to arrange all vertices along a circle in the order $\psi$ and join the edges of $G$ as chords. Let $S$ be a color set and $|S|=m$, a triple $(G,\psi,c)$ is an $m$-page \emph{book embedding} if $c:E(G) \rightarrow S$ is an edge-coloring such that $c(e^{\prime})\neq c(e^{\prime\prime})$ when $e^{\prime}$ and $e^{\prime\prime}$ cross in $\Psi$.~The \emph{book thickness}~$bt(G)$ of $G$ is the minimum integer $m$ such that an $m$-page book embedding exists. 
A book embedding $(G,\psi,c)$ is \emph{matching} if the edge-coloring $c$ is proper. The \emph{matching book thickness} $mbt(G)$ of $G$ is the minimum integer $m$ such that an $m$-page matching book embedding exists. A graph $G$ is called \emph{dispersable} if $mbt(G)=\Delta(G)$, \emph{nearly dispersable} if $mbt(G)=\Delta(G)+1$.

The $s\times t$ toroidal grid $T_{s,t}$ is the Cartesian product $C_{s}\Box C_{t}$ of two cycles $C_{s}, C_{t}$. If $s,t$ are both odd, we call $T_{s,t}$ an odd toroidal grid. The dispersability of the toroidal grid $T_{s,t}$ has aroused widespread interest.~In 2011,~Kainen ${[2]}$ proved that it is dispersable when $s,t$ are both even,~and nearly dispersable when $s$ is even and $t$ is odd.~Assume $s$ is odd and $s\geq t$, the case when $t = 3$ was solved by Shao et al $[3]$ in 2020, and the case when $t = 5$ was revealed by Joslin, Kainen and Overbay $[4]$ in 2021. Additionally, the case when $t$ is odd and $t\geq7$ was solved in $[5]$, which completely resolved the dispersability of the toroidal grid $T_{s,t}$.
In $[5]$, the authors divide two steps to find a 5-page matching book embedding of $T_{s,t}$. Firstly, the authors construct a 5-page matching book embedding $\Phi$ of $T_{t,t}$ such that some edges maintain a certain structure, and embed all edges recursively. The second step is to generalize $\Phi$ to a 5-page matching book embedding of $T_{s,t}$ by a trick of alternating orders based on the special structure of $\Phi$ in the first step.


In this work, we give a brief proof for the following main result in $[5]$.


\textbf{Theorem 1:} For all odd integers $s,t$ and $s \geq t \geq 7$, $mbt(T_{s,t})=5.$

\section{Proof of Theorem 1}




For brevity, we write $\mathbb{Z}_m$ for the set $\{1,2,\cdots,m\}$. Assume that $V(T_{s,t})=\{(p,q)~|~p\in \mathbb{Z}_s,q\in\mathbb{Z}_t\}$. For each $p\in \mathbb{Z}_s$, let the set $Y_p=\{(p,q)~|~q\in \mathbb{Z}_t\}$, which induces a $t$-cycle. Let $A_p=\{(p,1), (p,2), \cdots, (p,t)\}$, which is an arrangement of $Y_p$.

For $s$ is even, we construct a minimum-page matching book embedding $(T_{s,t},\psi,c)$ of $T_{s,t}$. Firstly, let $\psi=A_{1}^{-}A_{2}A_{3}^{-}A_{4}\cdots A_{s-1}^{-}A_{s}$; secondly, for each $i$ (see Fig.1), $i=2,4,\cdots,s$, the nested edges between $Y_{i-1}$ and $Y_i$ and between $Y_{i}$ and $Y_{i+1}$ are colored with red and blue, respectively; finally, by the edge chromatic number of cycles and the definition of $A_i$, the $t$-cycle induced by $Y_i$ can be colored with another two colors if $t$ is even, with another three colors if $t$ is odd, $i\in\mathbb{Z}_s$. 

\begin{figure}[htbp]
\centering
\includegraphics[width=0.6\textwidth]{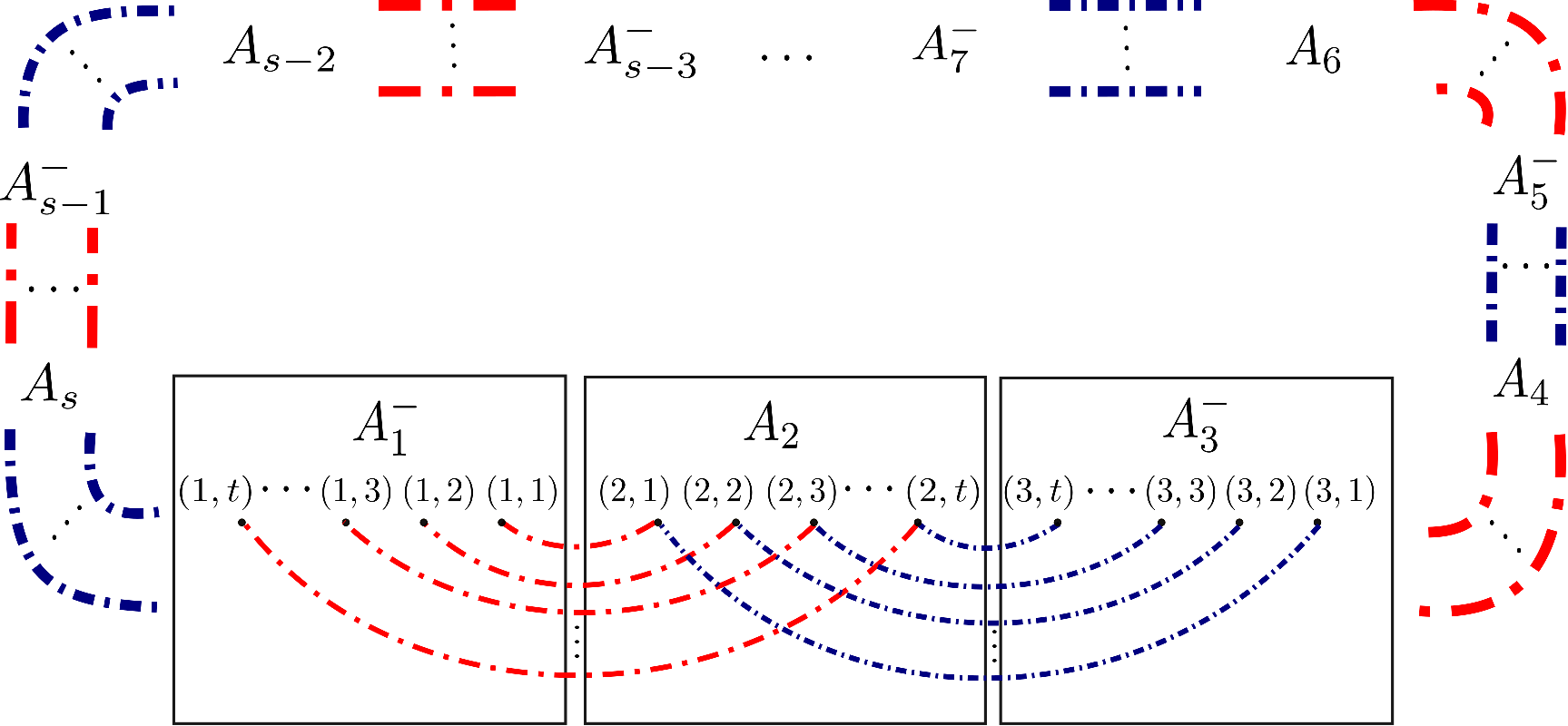}

Fig.1 ~The blue-page and red-page matching book embedding of $T_{n,t}$, $n$ is even.

\end{figure}
\vspace{-1em}

Next, we make use of the following two lemmas to give the proof of Theorem 1.

\noindent \textbf{Lemma 2.1.} $^{[2]}$
\emph{For all odd integers $s, t$,  $mbt(T_{s,t})\geq 5.$}



Our task is to construct a 5-page matching book embedding $\Psi$ of $T_{s,t}$. 
Inspired by the case where $s$ is even, our idea is to integrate the sets $Y_1,Y_2$ into a new set $U$, and to find an ordering of $U$ such that the vertex set of $T_{s,t}$ is divided into even disjoint sets $U,Y_3,Y_4\,\cdots,Y_{s-1},Y_s$. 


\noindent \textbf{Lemma 2.2. }\emph{For $s,t$ are both odd, $s \geq t \geq 7$, we have}
$mbt(T_{s,t})\leq 5.$

\noindent \textbf{Proof. }
Let $U$ be a rearrangement of $Y_{1}~\cup~ Y_{2}$. Specifically, $U=\{V_{2_1},V_{1_1},V_{1_2},V_{2_2},V_{2_3},V_{1_3},V_{1_4},V_{2_4},\\\cdots,V_{2_{t-2}},V_{1_{t-2}},V_{1_{t-1}},V_{2_{t-1}}, V_{2_{t}},V_{1_{t}}\}$, where $V_{i_j}$ is the $jth$ element of the ordered vertex set $A_i$, $i=1,2$. Put the vertices of $T_{s,t}$ counterclockwise along a circle in the order: $UA_3^{-}A_4A_5^{-}A_6\cdots A_{s-1}A_{s}^{-}.$
The edges of $T_{s,t}$ can be matching book embedded in five pages in the following two steps~(see Fig.2 for a $5$-page matching book embedding of $T_{11,9}$ and denote $(i,j)$ by $ij$ for short):

\begin{figure}[htbp]
\centering
\includegraphics[height=12.6cm ,width=0.88\textwidth]{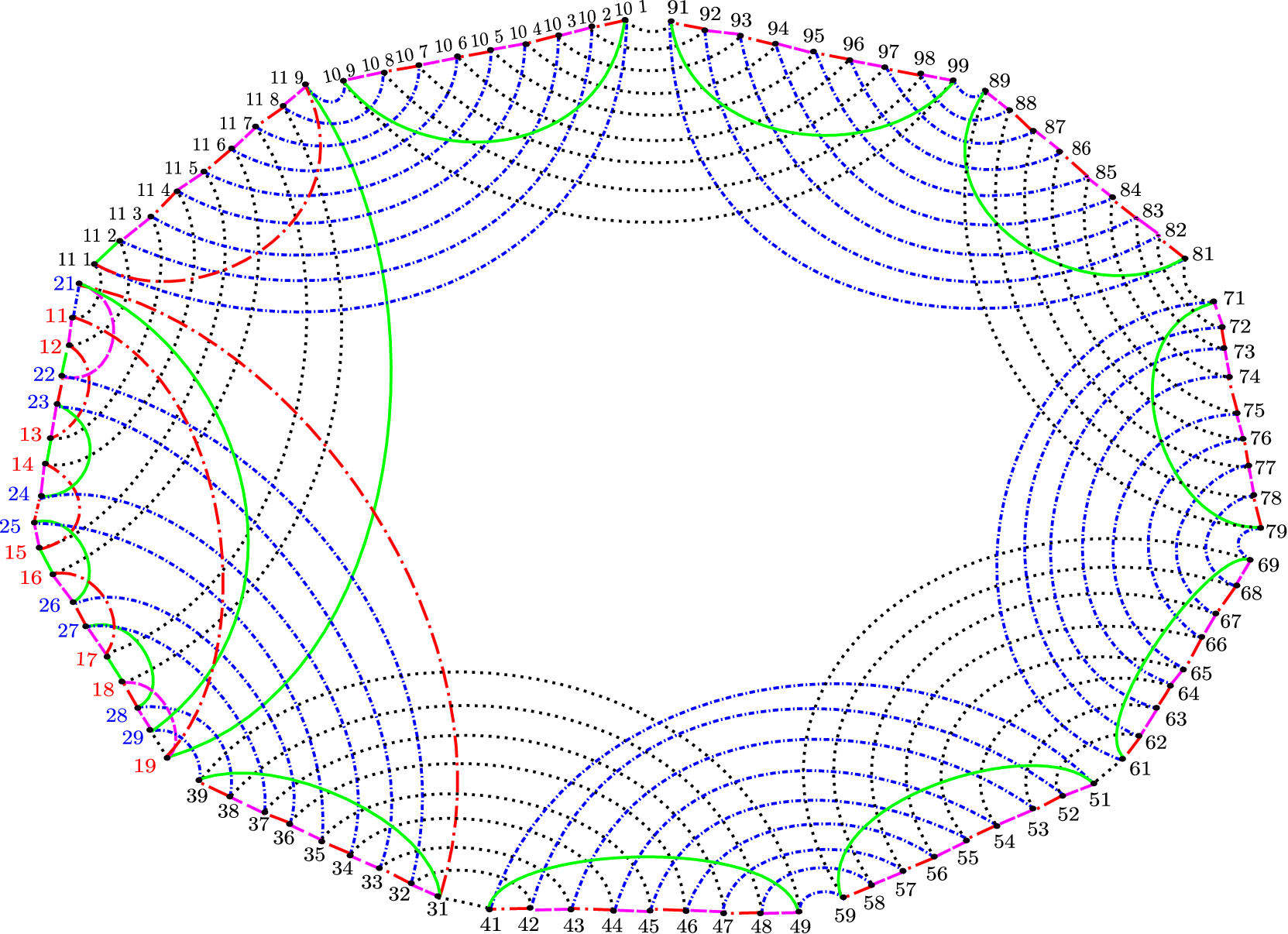}

Fig.2 ~A $5$-page matching book embedding of $T_{11,9}$.

\end{figure}


\noindent\textbf{Step 1:~}The coloring of the edges not belonging to the $t$-cycles induced by $Y_1,Y_2,\cdots,Y_s$.

Blue: $\{((i,j),(i+1,j))~|~i\in\{2,4,\cdots,s-1\},1\leq j\leq t,(i,j)\neq (2,1)\}\cup\{((1,1),(2,1))\}$;

Black: $\{((i,j),(i+1,j))~|~i\in\{3,5,\cdots,s\},1\leq j\leq t,(i,j)\neq (s,t)\}\cup\{((1,t),(2,t))\}$;

Green: $\{((1,t),(s,t))\}\cup\{((1,2),(2,2))\}$;

Red: $\{((2,1),(3,1))\}\cup\{((1,t-1),(2,t-1))\}$;

Purple: $\{((1,j),(2,j))~|~3\leq j\leq t-2\}$.

\noindent\textbf{Step 2:~}The coloring of the edges belonging to the $t$-cycles induced by $Y_1,Y_2,\cdots,Y_s$.

For the $t$-cycles induced by $Y_1$ and $Y_2$, the edges of $\{((i,j),(i,j+1))~|~i\in\{1,2\},j\in\{1,t-1\}\}$ are colored with purple, the edges of $\{((i,j),(i,j+1))~|~i\in\{1,2\},j\in\{2,4,\cdots,t-3\}\}\cup \{((1,1),(1,t))\}$ with red, and the remaining edges with green. For the $t$-cycles induced by $Y_3,Y_4,\cdots,Y_s$, it is easy to use red, green and purple to color these odd cycles well.

In summary, the result is established.
\hfill{$\square$}





%
%

\section*{Acknowledgment}

\indent


This work was partially funded by Science Research Project of Hebei Education Department,~China~(No.~ZD2020130) and the Natural Science Foundation of Hebei Province,~China~(No. A2021202013).

\end{document}